\documentclass[11pt, reqno]{amsart}
\usepackage{amsmath, amsthm, amscd, amsfonts, amssymb, graphicx, color}
\usepackage[bookmarksnumbered, colorlinks, plainpages]{hyperref}

\makeatletter \oddsidemargin.9375in \evensidemargin \oddsidemargin
\newtheorem{theorem}{Theorem}[section]
\newtheorem{lemma}[theorem]{Lemma}
\newtheorem{proposition}[theorem]{Proposition}
\newtheorem{corollary}[theorem]{Corollary}
\theoremstyle{definition}

\newtheorem{remark}[theorem]{Remark}

\numberwithin{equation}{section}

\begin{document}
\setcounter{page}{1}

%%%%%%%%%%%%%%%%%%%%%%%%%%%%%%%%%%%%%%%%%%%%%%%
%% Please do not remove the following statement.
%%%%%%%%%%%%%%%%%%%%%%%%%%%%%%%%%%%%%%%%%%%%%%%
%\noindent \textbf{{\footnotesize By submitting this extended
%abstract to IMC44 I \textcolor[rgb]{1.00,0.00,0.00}{confirm }that
%(i) I and any other coauthor(s) are responsible for its content and
%its originality; (ii) any possible coauthors agreed to its submission to IMC44.}}\\[1.00in]
%%%%%%%%%%%%%%%%%%%%%%%%%%%%%%%%%%%%%%%%%%%%%%%
%% ==>start from here
%% Start from here
%%%%%%%%%%%%%%%%%%%%%%%%%%%%%%%%%%%%%%%%%%%%%%%

%%%%%%%%%%%%%%%%%%%%%%%%%%%%%%%%%%%%%%%%%%%%%%%%%%%%%%%%%%%%%%%%%%%%%
% Insert title of your article. Note: \title[short title]{full title}
%%%%%%%%%%%%%%%%%%%%%%%%%%%%%%%%%%%%%%%%%%%%%%%%%%%%%%%%%%%%%%%%%%%%%
\title{The BSE-property for vector-valued $L^p-$algebras}
%%%%%%%%%%%%%%%%%%%%%%%%%%%%%%%%%%%%%%
% Author's name must be inserted here
%%%%%%%%%%%%%%%%%%%%%%%%%%%%%%%%%%%%%%
\author[Fatemeh Abtahi, Mitra Amiri and Ali Rejali]{Fatemeh Abtahi, Mitra Amiri and Ali Rejali}

%%%%%%%%%%%%%%%%%%%%%%%%
% Addresses
%%%%%%%%%%%%%%%%%%%%%%%%
%\authorsaddresses{ Department of Mathematics, University of Esfahan}
%%%%%%%%%%%%%%%%%%%%%%%%%%%%
% Subject class; see http://www.ams.org/mathscinet/msc/msc2010.html
%%%%%%%%%%%%%%%%%%%%%%%%%%%%,

%%%%%%%%%%%%%%%%%%%%%%%%%%%%%%%%%%%%%%%%%%%%%%%%%%%%%%%%%%%%%%%%%%%%%%%%%%%%%%%%%%%%%%%%%%
% keywords; Note that the number of keywords must be at least 3 items and at most 5 items.
%%%%%%%%%%%%%%%%%%%%%%%%%%%%%%%%%%%%%%%%%%%%%%%%%%%%%%%%%%%%%%%%%%%%%%%%%%%%%%%%%%%%%%%%%%

\begin{abstract}
Let $\mathcal A$ be a separable Banach algebra, $G$ be a locally compact Hausdorff group
and $1< p<\infty$. In this paper, we first provide a necessary and sufficient condition, for which
$L^p(G,\mathcal A)$ is a Banach algebra, under convolution product. Then we
characterize the character space of $L^p(G,\mathcal A)$, in the case where $\mathcal A$ is commutative
and $G$ is abelian. Moreover, we investigate the BSE-property for $L^p(G,\mathcal A)$ and
prove that $L^p(G,\mathcal A)$ is a BSE-algebra if and only if $\mathcal A$ is a BSE-algebra and $G$ is finite. Finally, we study the BSE-norm property for $L^p(G,\mathcal A)$ and show that if $L^p(G,\mathcal A)$ is a BSE-norm algebra then $\mathcal A$ is so. We prove the converse of this statement for the case where
$G$ is finite and $\mathcal A$ is unital.
\end{abstract}

\maketitle

\section{Introduction and preliminaries}

Let $(\mathcal A,\|\cdot\|_{\mathcal A})$ be a without order
commutative Banach algebra, in the sense that $a{\mathcal A}=\{0\}$
implies $a=0$ $(a\in {\mathcal A})$. We denote by
$\Delta({\mathcal A})$, the space consisting of all nonzero multiplicative linear
functionals on $\mathcal A$, called the character space of $\mathcal A$.
Throughout the paper, we assume that $\Delta({\mathcal A})$ is nonempty.
It shoulde be noted that $\Delta(\mathcal A)$, equipped with the weak$^*$topology,
inherited from ${\mathcal A}^*$, is a locally compact Hausdorff space.
We denote by $C_b(\Delta(\mathcal A))$ the space consisting of all bounded and continuous
functions on $\Delta(\mathcal A)$.
A bounded net $\{e_{\alpha}\}_{\alpha\in I}$ in $\mathcal A$, is called a bounded
$\Delta$-weak approximate identity for $\mathcal A$ if
$\lim_\alpha\varphi(ae_\alpha)=\varphi(a),$ for all $a\in\mathcal A$ and $\varphi\in\Delta(\mathcal A)$;
see \cite{JL}. Following \cite{Larsen}, a bounded linear operator
$T$ on $\mathcal A$, satisfying $T(ab)=aT(b)$, $(a,b\in \mathcal A)$,
is called a multiplier. The set of all multipliers on $\mathcal A$ is
denoted by $M({\mathcal A})$, which is a unital commutative Banach
algebra, called the multiplier algebra of ${\mathcal A}$.
By \cite[Theorem 1.2.2]{Larsen}, for any $T\in M({\mathcal A})$,
there exists a unique function $\widehat{T}\in C_b(\Delta(\mathcal A))$
such that $$\widehat{T(a)}(\varphi)=\widehat{T}(\varphi)\widehat{a}(\varphi),$$ for
all $a\in {\mathcal A}$ and $\varphi\in \Delta({\mathcal A})$.
A bounded and continuous function
$\sigma$ on $\Delta({\mathcal A})$ is called a BSE-function if
there exists a constant $\beta>0$ such that following the
inequality holds:
\begin{equation}\label{e15}
\left\vert \displaystyle\sum_{k=1}^n\alpha_k\sigma(\varphi_k)\right\vert\leq \beta\;\left\Vert
\displaystyle\sum_{k=1}^n\alpha_k\varphi_k\right\Vert_{{\mathcal
A}^*},
\end{equation}
for any finite number of  $\alpha_1,...,\alpha_n\in\Bbb C$ and the same number of
$\varphi_1,...,\varphi_n\in \Delta({\mathcal A})$.
The BSE-norm of $\sigma$ is the infimum of all $\beta$, which satisfying \eqref{e15}
and will be denoted by $\Vert\sigma\Vert_{BSE}$. By \cite{Tak1}, $C_{BSE}(\Delta({\mathcal A}))$,
the space consisting of all BSE-functions on $\Delta(\mathcal A)$, is a commutative and semisimple
Banach algebra, equipped with the norm $\Vert\cdot\Vert_{BSE}$ and pointwise product.
Then ${\mathcal A}$ is called a BSE-algebra (or has the BSE-property)
$\widehat{M({\mathcal A})}=C_{BSE}(\Delta({\mathcal A})),$ where
$\widehat{M({\mathcal A})}=\{\widehat{T}:T\in M({\mathcal
A})\}.$ By \cite[Corollary 5]{Tak1},
$\widehat{{M}({\mathcal A})}\subseteq  C_{BSE}(\Delta({\mathcal A})),$
if and only if ${\mathcal A}$ has a bounded
$\Delta$-weak approximate identity. It follows that all BSE-algebras possesses a bounded $\Delta$-weak
approximate identity. The Gelfand mapping $\Gamma_{\mathcal A}: \mathcal A\rightarrow C_b(\Delta(\mathcal A))$ is defined by
$a\mapsto\widehat{a},$ such that $\widehat{a}$ is the Gelfand
transform of $a$. The Banach algebra  $\mathcal A$
is called semisimple if $ker(\Gamma_{\mathcal A})=\{0\}.$
It should be noted that we always have $\widehat{\mathcal A}\subseteq
C_{BSE}(\Delta({\mathcal A}))$ and for any
$a\in\mathcal A$,
$$\|\widehat{a}\|_\infty\leq\|\widehat{a}\|_{BSE}\leq\|a\|_{\mathcal A},$$
where $\widehat{\mathcal A}$ is the range
of Gelfand mapping of $\mathcal A$. Let
${\mathcal M}({\mathcal A})$ is the normed algebra, consisting of all $\Phi\in C_b(\Delta(\mathcal A))$ such that $\Phi\cdot\widehat{\mathcal A}\subseteq\widehat{\mathcal A}.$
If $\mathcal A$ is semisimple then $\widehat{M({\mathcal
A})}={\mathcal M}({\mathcal A})$ \cite[p. 30]{Larsen}. Consequently, a semisimple and
commutative Banach algebra $\mathcal A$ is a BSE-algebra if
$C_{BSE}(\Delta({\mathcal A}))=\mathcal{M}({\mathcal A}).$

BSE-algebras and BSE functions were introduced
and investigated by Takahashi and Hatori \cite{Tak1}, and then
by some other authors, for various kinds of Banach algebras.
The acronym "BSE" stands for Bochner-Schoenberg-Eberlein and refers to a famous
theorem, for the additive group of real numbers, proved by
Bochner and Schoenberg \cite{Boch,SHo}. Then the result was generalized by Eberlein
\cite{Ebe}, for an abelian locally compact group $G$, which is indicating the
BSE-property of the group algebra $L^1(G)$ \cite{Rud}.
This result, has led Takahashi and Hatori \cite{Tak1} to introduce the BSE-property for any
commutative and without order Banach algebra $\mathcal A$. The interested reader is referred to \cite{AKT},
\cite{Dales}, \cite{FN}, \cite{IT}, \cite{IT1}, \cite{I}, \cite{Kan}, \cite{K},
\cite{Tak2}, and \cite{Tak3}.

In this paper, we investigate the BSE-property for $L^p(G,\mathcal A)$.
In fact, we first provide a necessary and sufficient condition, under which
$L^p(G,\mathcal A)$ is a Banach algebra, with convolution product. In fact, we generalize
$L^p$-conjecture \cite{S}, for the vector-valued case. We also
characterize the character space of $L^p(G,\mathcal A)$, in the case where $\mathcal A$ is commutative
and $G$ is abelian. Moreover, we prove that $L^p(G,\mathcal A)$ is a BSE-algebra
if and only if $\mathcal A$ is a BSE-algebra and $G$ is finite.
 Finally, we study the BSE-norm property for $L^p(G,\mathcal A)$ and show that if $L^p(G,\mathcal A)$
is a BSE-norm algebra then $\mathcal A$ is so. We prove the converse of this statement for the case where
$G$ is finite and $\mathcal A$ is unital.

\section{\bf Some basic results}

Let $(\mathcal A,\|\cdot\|_{\mathcal A}) $ be a Banach algebra,  $G$ be a locally compact Hausdorff group with a left Haar measure $\lambda$ and if $1\leq p<\infty $. A function $f: G\longrightarrow\mathcal A$ is called strongly measurable, if $f$ is Borel measurable and also $f(G)$ is separable in $\mathcal A$. Thus in the case where $\mathcal A$ is separable, the concept of measurability and strong measurability are equivalent. Most of the properties of integral theory in the complex version, are also valid for the vector-valued case. We refer to \cite{D}, as a complete survey in this issue.

Throughout the paper, $\mathcal A $ is a Banach algebra,  $G$ is a locally compact
Hausdorff group with a left Haar measure $\lambda$ and $1\leq p<\infty$. In the case where
$1<p<\infty$, we assume in addition that $\mathcal A$ is separable.
For any Borel measurable function $f: G\rightarrow\mathcal A$, let
$$\lVert f \rVert_{p,\mathcal A}=\left(\int_{G} \lVert f(x) \rVert_{\mathcal A}^{p} d\lambda(x )\right)^{\frac{1}{p}}.$$
Then $L^p (G,\mathcal A)$ is the Banach space, consisting of all Borel measurable functions $f: G \longrightarrow \mathcal A$ such that $\|f \Vert_{p,\mathcal A}<\infty$. For each $ f \in L^p (G,\mathcal A) $, define $ \overline{f}(x)= \Vert f(x) \Vert_{\mathcal A} $, for all $ x \in G $. Then $ f \in L^p (G,\mathcal A) $ if and only if $ \overline{f} \in L^p (G) $.
In this case, $\|\overline{f}\|_{p}=\|f\|_{p,\mathcal A}$. Recall that
for measurable vector-valued functions $f,g: G\longrightarrow\mathcal A$, the convolution multiplication
$$
(f *g)(x) =\int_Gf(y) g(y^{-1}x) d\lambda(y),
$$
is defined at each point $x \in G$ for which this makes sense.

In this section, we study some basic results about $L^p (G,\mathcal A)$. We give a necessary and sufficient condition for that $L^p (G,\mathcal A)$ is a
Banach algebra under convolution product. Note that by a famous conjecture, proved by Saeki \cite{S}, $L^p(G)$ is closed under convolution product if and only if $G$ is compact.

\begin{theorem}\label{t56}
Let $ \mathcal A $ be a separable Banach algebra, $G$ be a locally compact Hausdorff group and $1<p<\infty$. Then
$L^p (G,\mathcal A)$ is a Banach algebra under convolution product if and only if $ G $ is compact.
\end{theorem}

\begin{proof}
Let $L^p (G,\mathcal A)$ be a Banach algebra. To prove that $G$ is
compact, by \cite{S} it is sufficient to show that $ L^p (G) $ is closed under convolution. For any $ f,g \in L^p (G) $,
 define $\widetilde{f}(x)=f(x) a$ and $\widetilde{g}(x)=g(x) a$, where $a \in \mathcal A$ with $\Vert a \Vert =1$. Thus
 \begin{align*}
\left\Vert \widetilde{f} \ast \widetilde{g}\right\Vert_{p,\mathcal A}^{p} &=\int_{G}\left\Vert ( \widetilde{f} \ast \widetilde{g}) (x)\right\Vert_{\mathcal A}^{p} d \lambda (x)\\
&= \int_{G} \left\Vert  \int_{G} \widetilde{f} (y) \widetilde{g} ( y^{-1} x ) d \lambda (y)\right\Vert_{\mathcal A}^{p} d \lambda (x)\\
&=\int_{G} \left\Vert \int_{G} f(y) g( y^{-1} x )a^{2} d\lambda (y) \right\Vert_{\mathcal A}^{p} d \lambda (x)\\
&=\int_{G} \left\vert  ( f \ast g ) (x) \right\vert^{p} \left\Vert a^2\right\Vert_{\mathcal A}^{p} d \lambda (x)\\
&=\Vert a^2 \Vert_{\mathcal A}^{p} \,\, \Vert f \ast g \Vert_{p}^{p}.
\end{align*}
Consequently,
\begin{equation}\label{e1000}
\Vert\widetilde{f} \ast \widetilde{g} \Vert_{p,\mathcal A} = \Vert a^2 \Vert_{\mathcal A} \,\, \Vert  f \ast g  \Vert_{p}.
\end{equation}

By the hypothesis, we have
\begin{center}
$\Vert\widetilde{f}\ast \widetilde{g}\Vert_{p,\mathcal A} \leq\Vert \widetilde{f} \Vert_{p,\mathcal A } \, \Vert \widetilde{g} \Vert_{p,\mathcal A}$.
\end{center}
Thus by the equality \eqref{e1000} we obtain
$$
\Vert a^2 \Vert_{\mathcal A} \, \Vert f \ast g \Vert_{p}\leq\Vert a \Vert_{\mathcal A} \, \Vert f \Vert_{p} \, \Vert a \Vert_{\mathcal A} \, \Vert g \Vert_{p}=\Vert f \Vert_{p} \, \Vert g \Vert_{p}
$$
and so
\begin{center}
$ \Vert f \ast g \Vert_{p} \leq \frac{1}{\Vert a^2 \Vert_{\mathcal A}} \Vert f \Vert_{p}\, \Vert g \Vert_{p}.$
\end{center}
It follows that $L^p(G)$ is closed under convolution, which implies that $G$ is compact.

Conversely, suppose that $G$ is compact. Thus $L^p (G)$ is a Banach algebra under convolution product.
For all $f,g\in L^p(G,\mathcal A)$, we have  $\overline{f}, \overline{g}\in L^{p}(G) $ and
$$
\Vert \overline{f} \ast \overline{g} \Vert_{p}^{p} = \int_{G} \left\vert \int_{G} \Vert f(y) \Vert_{A} \,\, \Vert g ( y^{-1} x ) \Vert_{\mathcal A} d \lambda ( y )\right \vert^{p} d \lambda (x).
$$
Moreover, by a property of Bochner integrals we have
\begin{align*}
\Vert f \ast g \Vert_{p,\mathcal A}^{p}&=  \int_{G} \left\Vert \int_{G} f (y)g(y^{-1}x) d \lambda (y)\right\Vert_{\mathcal A}^{p} d \lambda (x)\\
&\leq \int_{G}\left(\int_{G} \Vert f(y) \Vert_{\mathcal A} \,\, \Vert g ( y^{-1} x ) \Vert_{\mathcal A} d \lambda (y)\right)^{p}d \lambda (x)\\
&= \Vert \overline{f} \ast \overline{g} \Vert_{p}^{p}.
\end{align*}
Consequently, we obtain
\begin{center}
$  \Vert f \ast g \Vert_{p,\mathcal A} \leq \Vert \overline{f} \ast \overline{g}  \Vert_{p} \leq \Vert \overline{f} \Vert_{p} \,\, \Vert \overline{g}  \Vert_{p}=\Vert f \Vert_{p,\mathcal A}\Vert g \Vert_{p,\mathcal A}.$
\end{center}
Therefore $ L^p (G,\mathcal A) $ is a Banach algebra, under convolution product.
\end{proof}

The following corollary is obtained by Theorem \ref{t56} and \cite{S}, immediately.

\begin{corollary}
Let $ \mathcal A $ be a separable Banach algebra, $G$ be a locally compact Hausdorff group and $1<p<\infty$.
Then the following statements are equivalent.
\begin{enumerate}
\item[(i)] $ L^p (G,\mathcal A)$ is a Banach algebra, under convolution product.
\item[(ii)] $L^p (G)$ is a Banach algebra, under convolution product.
\item[(iii)] $ G $ is compact.
\end{enumerate}
\end{corollary}

\begin{remark}\rm\label{r1}
Let $ \mathcal A $ be a separable Banach algebra, $G$ be a compact Hausdorff group and $1\leq p<\infty$.
By the Holder's inequality, we have
$ L^p(G)\subseteq L^1(G)$ and
$$
\|f\|_1=\int_G|f(x)|d\lambda(x)\leq\left(\int_G|f(x)|^pd\lambda(x)\right)^{1/p}=\|f\|_p,
$$
for all $f\in L^p(G)$. Now suppose that $\mathcal A$ is a separable Banach algebra and $f\in L^p(G,\mathcal A)$.
Then
$$
\|f\|_{1,\mathcal A}=\int_G||f(x)||d\lambda(x)\leq\left(\int_G||f(x)||^pd\lambda(x)\right)^{1/p}=\|f\|_{p,\mathcal A}.
$$
It follows that $f\in L^1(G,\mathcal A)$ and $\|f\|_{1,\mathcal A}\leq\|f\|_{p,\mathcal A}$.
In fact, we obtain $L^p(G,\mathcal A)\subseteq L^1(G,\mathcal A)$.
\end{remark}

It is known that $C_{00}(G)$, the space consisting of all complex-valued continuous functions on $G$ with compact support, is dense in $L^p(G)$ $(1\leq p<\infty)$. Now let
$$
C_{00}(G,\mathcal A)=\{f:G\longrightarrow\mathcal A,\;\;\text{f is continuous with compact support}\}.
$$
We denote it by $C(G,\mathcal A)$, whenever $G$ is compact.
In the next result, we investigate the same statement for the vector-valued case.
First, we introduce some notations, which we require for convenience. For any
$f\in L^p(G)$ and $a\in\mathcal A$, we denote by $f\otimes a$ the function, defined as
$$
f\otimes a(x)=f(x)a\;\;\;\;\;\;\;\;\;\;\;\;(x\in G).
$$
Then $f\in L^p(G,\mathcal A)$ and
$$\|f\otimes a\|_{p,\mathcal A}=\|f\|_p\|a\|_{\mathcal A}.$$
It should be noted that if $f\in C(G)$, then $f\otimes a$ belongs to $C(G,\mathcal A)$.
Let $C(G)\otimes\mathcal A$ be the vector space, generating by all $f\otimes a$, where $f\in C(G)$ and $a\in\mathcal A$. It is easily verified that
\begin{equation}\label{e1}
C(G)\otimes\mathcal A\subseteq C(G,\mathcal A)\subseteq L^p(G,\mathcal{A})\subseteq L^1(G,\mathcal{A}).
\end{equation}

\begin{proposition}\label{p1}
Let $ \mathcal A $ be a separable Banach algebra and $G$ be
a compact Hausdorff group. Then $C(G)\otimes\mathcal A$ is dense in $L^1(G,\mathcal{A})$.
\end{proposition}

\begin{proof}
Take $f \in L^1(G,\mathcal{A} )$. By \cite[Proposition 1.5.4.]{K},
there exist the sequences $\{f_n\}_{n}$ in $L^1(G)$ and also $\{a_{n}\}_{n}$ in $\mathcal{A}$ such that
$f= \sum_{n=1}^{\infty} f_n \otimes a_n $ and $\sum_{n=1}^{\infty}
\lVert f_{n}\rVert_{1} \lVert a_{n}\rVert_{\mathcal{A}}< \infty$. Since $C(G)$ is dense $L^1(G) $, for each $n\in\Bbb N$ and $\varepsilon>0$, there exists $g_{n}\in C(G)$ such that
$$\lVert g_{n}-f_{n}\rVert_{1} <\frac{\varepsilon}
{2^{n+1}(\lVert a_{n}\rVert_{\mathcal A}+1)}.$$
Define $g=\sum_{n=1}^{\infty} g_n \otimes a_n $. Then
\begin{eqnarray*}
\sum_{n=1}^{\infty}\left\Vert g_{n}\right\Vert_{1}\left\Vert a_{n}\right\Vert_{\mathcal{A}}&\leq&
\sum_{n=1}^{\infty}\left\Vert f_n-g_{n}\right\Vert_{1}\left\Vert a_{n}\right\Vert_{\mathcal{A}}+\sum_{n=1}^{\infty}\left\Vert f_{n}\right\Vert_{1} \left\Vert a_{n}\right\Vert_{\mathcal{A}}\\
&\leq&\sum_{n=1}^{\infty}\frac{\varepsilon\left\Vert a_{n}\right\Vert_{\mathcal{A}}}
{2^{n+1}(\left\Vert a_{n}\right\Vert_{\mathcal A} +1)}+\sum_{n=1}^{\infty}\left\Vert f_{n}\right\Vert_{1} \left\Vert a_{n}\right\Vert_{\mathcal{A}}\\
&\leq&\varepsilon+\sum_{n=1}^{\infty}\left\Vert f_{n}\right\Vert_{1} \left\Vert a_{n}\right\Vert_{\mathcal{A}}\\
&<&\infty,
\end{eqnarray*}
which implies that $g\in L^1(G,\mathcal A)$. Moreover,
\begin{align*}
\lVert g-f\rVert_{1,\mathcal A}&=\left\lVert \sum_{n=1}^{\infty}(g_{n}-f_{n})\otimes a_n\right\Vert_{1,\mathcal A}\\
&\leq \sum_{n=1}^{\infty}\lVert g_{n}-f_{n}\rVert_{1}\;\lVert a_n\rVert_{\mathcal A}\\
&\leq \sum_{n=1}^{\infty}\frac{\varepsilon\lVert a_{n}\rVert_{\mathcal{A}}}
{2^{n+1}(\lVert a_{n}\rVert_{\mathcal{A}}+1)}\\
&<\varepsilon.
\end{align*}
Now let $h_{n}=\sum_{k=1}^{n} g_k \otimes a_k$ $(n\in\Bbb N)$. Then
$\{h_{n}\}\subseteq C(G)\otimes\mathcal A$ and
$$
\lim_{n}\lVert h_{n}-g\rVert_{1,\mathcal{A}}=0.
$$
It follows that there exists $n\in\Bbb N$ such that
$\lVert h_{n}-g\rVert_{1,\mathcal{A}}<\varepsilon$. Hence, we obtain
$$\lVert f-h_{n} \rVert_{1,\mathcal{A}}\leq \lVert f-g\rVert_{1,\mathcal{A}}+\lVert g-h_{n}\rVert_{1,\mathcal{A}}<2\varepsilon.$$
This completes the proof.
\end{proof}

\begin{proposition}
Let $ \mathcal A $ be a separable commutative Banach algebra, $G$ be
an abelian compact Hausdorff group and $1\leq p<\infty$. Then $L^p(G,\mathcal A)$ is a dense ideal
in $L^1(G,\mathcal A)$.
\end{proposition}

\begin{proof}
Proposition \ref{p1} together with the inclusions \eqref{e1} imply that $L^p(G,\mathcal A)$ is dense in $L^1(G,\mathcal A)$.
We show that $L^p(G,\mathcal A)$ is an ideal in $L^1(G,\mathcal A)$. To this end, let $f\in L^1(G,\mathcal A)$ and $g\in L^p(G,\mathcal A)$. We indicate that $f*g\in L^p(G,\mathcal A)$. Since $\overline{f}\in L^1(G)$ and $\overline{g}\in L^p(G)$, together with the fact that
$L^p(G)$ is an ideal in $L^1(G)$, we obtain that $\overline{f}*\overline{g}\in L^p(G)$.
Moreover,
$$
\|f*g\|_{p,\mathcal A}\leq\|\overline{f}*\overline{g}\|_p<\infty.
$$
It follows that $f*g\in L^p(G,\mathcal A)$, as claimed.
\end{proof}

The next result is obtained from \cite[Theorem 2.3]{Ba}. We refer to \cite{Ba},
for the general definition of abstract Segal algebras.

\begin{corollary}\label{c1}
Let $\mathcal A$ be a separable commutative Banach algebra, $G$ be
an abelian compact Hausdorff group and $1\leq p<\infty$. Then $L^p(G,\mathcal A)$ is an
abstract Segal algebra in $L^1(G,\mathcal A)$.
\end{corollary}

\begin{remark}
Let $\mathcal A$ be a separable commutative Banach algebra with a bounded approximate identity, $G$ be
an abelian compact Hausdorff group and $1< p<\infty$. Then by \cite[Theorem 2.9.21]{Dales1}
and \cite[Proposition 1.5.4]{K}, $L^1 (G,\mathcal A)$ possesses a bounded approximate
identity. Now Cohen factorization theorem implies that
$L^1(G,\mathcal A)*L^p(G,\mathcal A)=L^p(G,\mathcal A)$.
\end{remark}

Let $\mathcal A$ be a commutative Banach algebra and $G$ be an abelian locally compact Hausdorff group
with dual group $\widehat{G}$, introduced in \cite{Rud}. By \cite[Proposition 1.5.4]{K}, $L^1(G,\mathcal A)$ is
isomorphic with the projective tensor product $L^1(G)$ and $\mathcal A$.
Moreover, by \cite[Theorem 2.11.2]{K}, the character space of $L^1(G,\mathcal A)$ is
homeomorphic with $\widehat{G}\times\Delta(\mathcal A)$, such that for all $\chi\in\widehat{G}$ and $\varphi\in\Delta(\mathcal A)$,
$$
\chi\otimes\varphi:L^1(G)\widehat{\otimes}\mathcal A\longrightarrow\Bbb C,
$$
defined as, $$(\chi\otimes\varphi)(f\otimes a)=\chi(f)\varphi(a)\;\;\;\;(f\in L^1(G),a\in\mathcal A).$$
Moreover, for all $f\in L^1(G,\mathcal A)$, we have
$$
(\chi\otimes\varphi)(f)=\varphi\left(\int_G\overline{\chi(x)}f(x)d\lambda(x)\right)=
\int_G\overline{\chi(x)}\varphi(f(x))d\lambda(x).
$$
Now Corollary \ref{c1} together with \cite[Lemma 2.2]{ANR} yield the next result.

\begin{proposition}
Let $\mathcal A$ be a separable commutative Banach algebra, $G$ be
an abelian compact Hausdorff group and $1< p<\infty$. Then
$$
\Delta(L^p(G,\mathcal A))=\{(\chi\otimes\varphi)|_{L^p(G,\mathcal A)}:\;\;\chi\in\widehat{G},\varphi\in\Delta(\mathcal A)\}.
$$
\end{proposition}

\section{\bf $L^p(G,\mathcal A)$ as a BSE-algebra }

Let $G$ be an abelian compact Hausdorff group and $1<p<\infty$. Then $L^p(G)$
is a commutative Banach
algebra, which is an ideal in its second dual. Now by \cite[Theorem 3.1]{Kan}, the following assertions are equivalent.
\begin{enumerate}
\item[(i)] $L^p(G)$ is a BSE-algebra.
\item[(ii)] $L^p(G)$ possesses a bounded $\Delta-$weak approximate identity.
\item[(iii)] $L^p(G)$ admits a bounded approximate identity.
\end{enumerate}
Since $L^p(G)$ is reflexive, the above equivalent conditions imply that $L^p(G)$ is unital and so $G$ is finite.

In this section, we investigate the BSE-property for $L^p(G,\mathcal A)$.
We commence with the following proposition.

\begin{proposition}
Let $\mathcal A$ be a separable commutative Banach algebra,
$G$ be an abelian compact Hausdorff group and $1<p<\infty$.
Then $L^p (G,\mathcal A)$ is semisimple if and only if $\mathcal A $ is semisimple.
\end{proposition}

\begin{proof}
Let $L^p(G,\mathcal A)$ be semisimple. Then $\Delta (L^p(G,\mathcal A))=
\widehat{G} \times \Delta(\mathcal A)$ separates the points of $L^p(G,\mathcal A)$.
Take $ a,b \in \mathcal A$ with $a \neq b$. Then
\begin{center}
$ 1_{G}\otimes a \neq  1_{G}\otimes b,$
\end{center}
where $1_G$ is the constant function $1$ on $G$. By the assumption,
there exist $\chi\in\widehat{G}$ and  $\varphi\in \Delta(\mathcal{A})$ such that
\begin{center}
$(\chi \otimes\varphi)(1_G\otimes a)\neq(\chi\otimes\varphi)(1\otimes b).$
\end{center}
Consequently,
\begin{center}
$\chi (1)\varphi(a)\neq \chi(1)\varphi(b),$
\end{center}
and so $\varphi(a)\neq\varphi(b) $. Thus $ \mathcal{A} $ is semisimple.

Conversely, suppose that $\mathcal A$ is semisimple and take $f,g \in L^p(G,\mathcal A)$ with $f\neq g$.
Then $\lVert f-g\rVert_{p,\mathcal A}\neq 0$. Thus $\lambda(\{x\in G: f(x)\neq g(x)\})\neq 0.$
It follows that $\lVert f-g\rVert_{1,\mathcal A}\neq 0.$
By the semisimplicity of
$L^1(G,\mathcal A)$ \cite[Theorem 2.11.8]{K}, there exists
$\psi \in \Delta(L^1(G, \mathcal A))=\Delta( L^p(G, \mathcal A)) $
such that $\psi(f)\neq \psi(g)$. It follows that $L^p(G,\mathcal A)$ is semisimple.
\end{proof}

The next proposition is applied in some further results.
The proof is straightforward and left to the reader.

\begin{proposition}\label{p20}
Let $\mathcal A$ be a separable commutative unital Banach algebra  with the identity element $e$ such that
$\|e\|_{\mathcal A}=1$, $G$ be
an abelian compact Hausdorff group and $1\leq p<\infty$.
Then for any $\chi_1,\cdots,\chi_n\in\widehat{G}$ and $\varphi_1,\cdots,\varphi_n\in\Delta(\mathcal A)$ and
the same number of $c_1,\cdots,c_n\in \Bbb C$, we have
$$
\left\Vert \sum_{i=1}^{n} c_i \chi_i\right\Vert_{L^p (G)^{\ast}}\leq\left\Vert \sum_{i=1}^{n} c_i ( \chi_i \otimes \varphi_i )\right\Vert_{L_p (G, \mathcal{A} )^{\ast}}
$$
and
$$
\left\Vert \sum_{i=1}^{n} c_i \varphi_i\right\Vert_{{\mathcal A}^{\ast}}\leq\left\Vert \sum_{i=1}^{n} c_i ( \chi_i \otimes \varphi_i )\right\Vert_{L_p (G, \mathcal{A} )^{\ast}}
$$
\end{proposition}

\begin{proposition}\label{p90}
Let $\mathcal A$ be a separable commutative Banach algebra, $G$ be
an abelian compact Hausdorff group and $1<p<\infty$.  Then $L^p(G,\mathcal{A})$ has a
bounded $\Delta$-weak approximate identity if and only if $\mathcal{A}$ has so.
\end{proposition}

\begin{proof}
First suppose that $\{f_{\alpha}\}_{\alpha}$ is a bounded $\Delta$-weak approximate identity for
$L^p (G,\mathcal{A})$, bounded by $M>0$. Take $\chi_0 \in \widehat{G}$ to be fixed and for each $\alpha$ let
$$
e_{\alpha}:= \chi_0 ( f_{\alpha} )=\int_G f_{\alpha} (x) \overline{\chi_0 (x)}d\lambda(x).
$$
Then for each $\varphi \in \Delta ( \mathcal{A})$, we have
$$
\varphi (e_{\alpha} )=\int_G \varphi ( f_{\alpha}  (x)) \overline{\chi_0 (x)}d\lambda(x)=\chi_0 (\varphi \circ f_{\alpha})
=\chi_0 \otimes \varphi (f_{\alpha})\longrightarrow_\alpha 1 .
$$
Moreover,
$$
\Vert e_{\alpha}\Vert_{\mathcal{A}}\leq\int_G\Vert f_{\alpha} (x) \Vert_{\mathcal{A}}\vert \overline{\chi_0 (x)}\vert d\lambda(x)\leq M\int_G\overline{\chi_0 (x)}d\lambda(x)\leq M\lambda (G).
$$
Therefore $\{e_{\alpha}\}_\alpha$ is a bounded $\Delta$-weak approximate identity for $\mathcal{A}$.

Conversely, let $\mathcal A$ has a bounded $\Delta$-weak approximate identity, denoted by $\{e_\alpha\}_{\alpha}$ and $\{f_\beta\}_{\beta}$ be the
 bounded approximate identity for $L^1(G)$. Then for all $\varphi\in\Delta(\mathcal A)$,
$ \varphi(e_\alpha) \longrightarrow 1$ and so for all $\chi\in\widehat{G}$, we have
\begin{center}
$(\chi \otimes\varphi)(f_\beta \otimes e_{\alpha})\longrightarrow 1$.
\end{center}
 This implies that $ (f_\beta \otimes e_{\alpha})_{\alpha,\beta} $ is a
  bounded $\Delta$-weak approximate identity for $L^p(G,\mathcal{A})$.
\end{proof}

We state here the main result of this section.

\begin{theorem}\label{t5}
Let $\mathcal A $ be a semisimple commutative separable and unital Banach algebra,
$G$ be an abelian compact Hausdorff group and $1<p<\infty$. Then $L^p(G,\mathcal A)$
is a BSE-algebra if and only if $\mathcal A$ is a BSE-algebra and $G$ is finite.
\end{theorem}

\begin{proof}
Let $ L^p (G,\mathcal A) $ be a BSE-algebra. By \cite[Corollary 5]{Tak1}, $L^p (G,\mathcal A)$ admits
a bounded $\Delta-$weak approximate identity $\lbrace f_{\alpha} \rbrace_{\alpha}$, bounded by $M$. It follows that
\begin{center}
$\widehat{f}_{\alpha}(\chi\otimes\varphi) \longrightarrow_\alpha 1,$
\end{center}
and so
\begin{equation}\label{e3}
\int_{G}\overline{\chi(x)}\varphi (f_{\alpha}(x))dx \longrightarrow_\alpha 1,
\end{equation}
for all $\chi\in\widehat{G}$ and $\varphi\in\Delta (\mathcal A)$. Now, take $\varphi_0\in\Delta (\mathcal A)$ to be fixed and define $ g_\alpha (x) = \varphi_0(f_\alpha (x)) $ $(x\in G)$. Thus
\begin{eqnarray*}
\left(\int \vert  g_\alpha(x)\vert^{p} dx\right)^{1/p}&=&\left(\int_{G} \vert \varphi_0 ( f_\alpha (x))\vert^{p} dx\right)^{1/p}\\
&\leq& \Vert \varphi_0 \Vert \left(\int_{G} \Vert f_\alpha (x) \Vert^{p} dx\right)^{1/p}\\
&\leq& M.
\end{eqnarray*}
Consequently, $\lbrace g_\alpha \rbrace_{\alpha}$ is a bounded net in $L^p (G)$.
Moreover, by \eqref{e3}, for all $\chi\in \widehat{G}$ we obtain
$$
\widehat{g_{\alpha}}(\chi)=\chi(g_{\alpha})=\int_{G}g_{\alpha}(x)\overline{\chi(x)}dx
=\int_{G}\overline{\chi(x)}\varphi_0(f_{\alpha}(x))
\longrightarrow_\alpha 1.
$$
Thus $\lbrace g_\alpha \rbrace_{\alpha}$ is a bounded $\Delta-$weak
approximate identity for $L^p (G)$, which implies that $G$ is finite.
In the sequel, we show that $\mathcal A$ is a BSE-algebra.
Proposition \ref{p90} implies that $\mathcal A$
has a bounded $\Delta$-weak approximate identity. Thus
\cite[Corollary 5]{Tak1} implies that
\begin{equation}\label{e100}
\mathcal{M}(\mathcal A) \subseteq  C_{BSE}(\Delta (\mathcal{A})).
\end{equation}
Now, we prove that the reverse of the inclusion \eqref{e100}. Take $\sigma\in C_{BSE}(\Delta (\mathcal{A}))$ and define the function $\overline{\sigma}:\widehat{G}\times \Delta (A)\longrightarrow \mathbb C$ as
$\overline{\sigma}(\chi\otimes\varphi):=\sigma(\varphi),$
for all $\chi\in \widehat{G} $ and $ \varphi\in \Delta (\mathcal{A}) $. We show
that $ \overline{\sigma}\in  C_{BSE}(\Delta (L^p (G, \mathcal{A})))$. Note that since $G$  is finite, thus $G=\widehat{G}$.
Now Proposition \ref{p20} implies that for all complex numbers
$c_1,\cdots,c_n$ and the same number of $\chi_{1}\otimes\varphi_{1},
\cdots,\chi_n\otimes\varphi_n$ of $\widehat{G}\times\Delta(\mathcal A)$,
$$
\left\vert \sum_{i=1}^{n} c_i\overline{\sigma}(\chi_{i}\otimes\varphi_{i})\right\vert
\leq\Vert \sigma\Vert_{BSE}\left\Vert \sum_{i=1}^{n} c_i\varphi_{i}\right\Vert_{\mathcal{A}^{*}}
\leq\Vert \sigma\Vert_{BSE}\left\Vert \sum_{i=1}^{n} c_i
(\chi_{i}\otimes\varphi_{i})\right\Vert_{L^p (G, \mathcal{A})^*}.
$$
Consequently, $\overline{\sigma}\in C_{BSE}(\Delta(L^p (G,\mathcal{A})))=\mathcal {M}(L^p (G,\mathcal{A}))$.
Hence,  for each $f \in C(G)$ and $a\in\mathcal A $,
 there exists $ g \in L^p (G, \mathcal{A}) $ such that
$\overline{\sigma}\;\widehat{(f\otimes a)}=\widehat{g}.$
It follows that
$$\sigma(\varphi)\;\chi(f)\varphi(a)=\chi(\varphi \circ g).
\quad \quad (\chi\in \widehat{G}, \varphi\in \Delta(\mathcal A))$$
Taking $\chi\in\widehat{G}$ and $f\in C(G)$ with $\chi(f)=1 $, we obtain
$$
\sigma(\varphi)\widehat{a}(\varphi)=\chi(\varphi \circ g)
=\varphi\left(\int_G g(x)\overline{\chi(x)}d\lambda(x)\right)
:=\varphi(b),
$$
where $b:=\int_G g(x) \overline{\chi(x)}d\lambda(x)$.
Consequently, $\sigma\widehat{a}=\widehat{b}$, which implies that
$\sigma\in{\mathcal M}(\mathcal A)$. It follows that $\mathcal A$ is a $BSE-$algebra.

Conversely, suppose that $\mathcal A$ is a BSE-algebra and $G$ is finite. Thus by
Proposition \ref{p90}, $L^p(G,\mathcal A)$ possesses a bounded $\Delta$-weak approximate identity and so
$$
\mathcal{M}(L^p(G,\mathcal A))\subseteq C_{BSE}(\Delta(L^p(G,\mathcal A))),
$$
by \cite[Corollary 5]{Tak1}. To prove the reverse of the inclusion, suppose that $\sigma\in C_{BSE}(\Delta(L^p(G,\mathcal A)))$
and $f \in L^p(G,\mathcal A)$. We should find $g \in L^p(G,\mathcal A)$
such that $\sigma\widehat{f}=\widehat{g}$. Note that since $G$ is finite, it follows that
$G=\widehat{G}$. For each $\chi\in\widehat{G}$, define
$$\sigma_{\chi}: \Delta(\mathcal A)\rightarrow\mathbb{C}$$
with
$$\sigma_{\chi}(\varphi)=\sigma(\chi\otimes\varphi)\;\;\;\;\;(\varphi\in\Delta(\mathcal A)).$$
It is easily verified that $\sigma_{\chi}\in C_b(\Delta(\mathcal A))$. Now we show that $\sigma_{\chi}\in C_{BSE}(\Delta(\mathcal A))$.
For any $c_1,...,c_n \in \mathbb{C}$ and the same number of $\varphi_1,...,\varphi_n$ in $\Delta(\mathcal A)$, we easily obtain
$$\left|\sum_{i=1}^{n}c_i\sigma_{\chi}(\varphi_i)\right|\leq\|\sigma\|_{BSE}\left\|\sum_{i=1}^{n}c_i\varphi_i \right\|_{{\mathcal A}^*}.$$
Consequently, $\sigma_{\chi}\in C_{BSE}(\Delta(\mathcal A))$ and since $\mathcal A$ is a BSE algebra, $\sigma_{\chi}\in{\mathcal M}(\mathcal A)$.
It follows that for each $\chi\in\widehat{G}$, there exists $a_{\chi}\in\mathcal A$ such that
$$
\sigma_{\chi}\widehat{f(\chi)}=\widehat{a_{\chi}}.
$$
Now define the function
$$g:G=\widehat{G}\rightarrow \mathcal A\;\;\;\;\;\;\;g(\chi):=a_{\chi},$$
which belongs to
$L^p(G,\mathcal A)$. Finally, for each $\chi\otimes\varphi$ in $\Delta(L^p(G,\mathcal A))$, we have
$$
\sigma(\chi\otimes\varphi)\widehat{f}(\chi\otimes\varphi)=\sigma_\chi(\varphi)\widehat{f(\chi)}(\varphi)
=\widehat{g(\chi)}(\varphi)
=\widehat{g}(\chi\otimes\varphi).
$$
Consequently, $\sigma\widehat{f}=\widehat{g}$. It follows that $\sigma\in{\mathcal M}(L^p(G,\mathcal A))$.
Therefore $L^p(G,\mathcal A)$ is a BSE algebra.
\end{proof}

\section{\bf $L^p(G,\mathcal A)$ as a BSE-norm algebra }

It is known that in any commutative Banach algebra
$\mathcal A$,
\begin{equation}\label{e13}
\|\widehat{a}\|_\infty\leq\|\widehat{a}\|_{BSE}\leq\|a\|_{\mathcal A},
\end{equation}
for all $a\in \mathcal A$. In \cite{Dales2}, the authors
investigated a class of commutative Banach algebras which satisfies
the following condition $$\|a\|_{\mathcal A}\leq K\|\widehat{a}\|_{BSE}\;\;\;\;\;\;\;\;(a\in\mathcal A),$$
for some $K>0$, called BSE-norm algebras. They indicates that $L^p(G)$ is always a BSE-norm algebra.
As the final results, we show that this result is not necessarily valid for the vector-valued case.
We commence with the following elementary lemma.

\begin{lemma}\label{l1}
Let $ \mathcal A $ be a separable and semisimple Banach algebra, $G$ be
a compact Hausdorff group and $1\leq p<\infty$. Then
$ L^p (G,\mathcal A) $ is unital if and only if $\mathcal A$ is unital and $G$ is finite.
\end{lemma}

\begin{proof}
Suppose that $L^p (G,\mathcal A)$ is unital. Thus $\widehat{G}\times\Delta(\mathcal A)$
is compact and so $\Delta(\mathcal A)$ and  $\widehat{G}$ are compact. It follows that $\mathcal A$ is unital and $G$ is discrete, which implies that $G$ finite. The converse is similar.
\end{proof}

\begin{theorem}\label{t57}
Let $\mathcal A $ be a semisimple commutative separable unital Banach algebra, $G$ be an abelian compact Hausdorff group and $1\leq p<\infty$. If $L^p(G,\mathcal A)$ is a BSE-norm algebra then $\mathcal A$ is so. The converse is true, whenever $G$ is finite and $\mathcal A$ is unital.
\end{theorem}	

\begin{proof}
By the hypothesis, there exists $M>0$ such that
$$
\|f\|_{p,\mathcal A}\leq M\|\widehat{f}\|_{BSE}\;\;\;\;\;\;(f\in L^p(G,\mathcal A)).
$$
Taking $a\in\mathcal A$, for the constant fuction $f_a$ we have
\begin{equation}\label{e5}
\|a\|_{\mathcal A}=\|f_a\|_{p,\mathcal A}\leq M\|\widehat{f_a}\|_{BSE}\;\;\;\;\;\;(f\in L^p(G,\mathcal A)).
\end{equation}
We show that $\|\widehat{f_a}\|_{BSE}\leq\|\widehat{a}\|_{BSE}$. For $\chi_1\otimes\varphi_1,\cdots,\chi_n\otimes\varphi_n\in\Delta(L^p(G,\mathcal A))$ and the same number of $c_1,\cdots,c_n\in\Bbb C$ we have
\begin{eqnarray}\label{e2}
\left|\sum_{i=1}^nc_i\widehat{f_a}(\chi_i\otimes\varphi_i)\right|&=&\left|\sum_{i=1}^nc_i(\chi_i\otimes\varphi_i)(f_a)\right|\nonumber\\
&=&\left|\sum_{i=1}^n\left(c_i\int_G\overline{\chi_i(x)}d\lambda(x)\right)\varphi_i(a)\right|\nonumber\\
&\leq&\|\widehat{a}\|_{BSE}\left\|\sum_{i=1}^n\left(c_i\int_G\overline{\chi_i(x)}d\lambda(x)\right)\varphi_i\right\|_{{\mathcal A}^*}.
\end{eqnarray}
On the other hand,
\begin{eqnarray}\label{e3}
\left\|\sum_{i=1}^nc_i(\chi_i\otimes\varphi_i)\right\|_{L^p(G,\mathcal A)^*}&=&
\sup_{\|f\|_{p,\mathcal A}\leq 1}\left|\sum_{i=1}^nc_i\int_G\overline{\chi_i(x)}\varphi_i(f(\chi_i))d\lambda(x)\right|\nonumber\\
&\geq&\sup_{\|a\|_{\mathcal A}\leq 1}\left|\sum_{i=1}^nc_i\int_G\overline{\chi_i(x)}\varphi_i(a)d\lambda(x)\right|\nonumber\\
&=&\left\|\sum_{i=1}^n\left(c_i\int_G\overline{\chi_i(x)}d\lambda(x)\right)\varphi_i\right\|_{{\mathcal A}^*}.
\end{eqnarray}
Now \eqref{e2} and \eqref{e3} imply that
$$
\left|\sum_{i=1}^nc_i\widehat{f_a}(\chi_i\otimes\varphi_i)\right|\leq\|\widehat{a}\|_{BSE}
\left\|\sum_{i=1}^nc_i(\chi_i\otimes\varphi_i)\right\|_{L^p(G,\mathcal A)^*}.
$$
It follows that
\begin{equation}\label{e4}
\|\widehat{f_a}\|_{BSE}\leq\|\widehat{a}\|_{BSE}.
\end{equation}
By \eqref{e5} and \eqref{e4} we obtain
$$
\|a\|_{\mathcal A}\leq M\|\widehat{a}\|_{BSE},
$$
which implies that $\mathcal A$ is a BSE-norm algebra.

Conversely, suppose that $G$ is finite and $\mathcal A$ is unital. By Lemma \ref{l1}, $L^p(G,\mathcal A)$ is unital.
Moreover, Theorem \ref{t5} implies that $L^p(G,\mathcal A)$ is a BSE-algebra. Thus $L^p(G,\mathcal A)$ is a BSE-norm algebra, by \cite[Theorem page 40]{Dales2}.
\end{proof}

\vspace{9mm}

{\footnotesize \noindent
 F. Abtahi\\
Department of Pure Mathematis\\
Faculty of Mathematics and Statistics\\
University of Isfahan\\
Isfahan, Iran\\
    f.abtahi@sci.ui.ac.ir\\
     abtahif2002@yahoo.com\\

\noindent
M. Amiri\\
Department of Pure Mathematis\\
Faculty of Mathematics and Statistics\\
University of Isfahan\\
Isfahan, Iran\\
m.amiri@sci.ui.ac.ir\\
mitra75amiri@gmail.com\\

\noindent
 A. Rejali\\
 Department of Pure Mathematis\\
Faculty of Mathematics and Statistics\\
University of Isfahan\\
Isfahan, Iran\\
    rejali@sci.ui.ac.ir\\
   a.rejali20201399@gmail.com

\end{document}